\documentclass{llncs}
\begin{document}

\pagenumbering{arabic}
\pagestyle{headings}
\def\sof{\hfill\rule{2mm}{2mm}}
\def\ls{\leq}
\def\gs{\geq}
\def\tx{{ \left( \frac1{2\sqrt{x}} \right)}}

\title{ {\sc Restricted set of patterns, continued fractions, and Chebyshev polynomials} }
   
\author{Toufik Mansour}
\institute{Department of Mathematics, \\ 
		University of Haifa, Haifa, Israel 31905\\
		{\tt tmansur11@hotmail.com }}
\maketitle
\section*{\centering{\sc Abstract}}
We study generating functions for the number of permutations in $S_n$
subject to set of restrictions. One of the restrictions belongs to $S_3$,
while the others to $S_k$. It turns out that in a large variety of cases
the answer can be expressed via continued fractions, and Chebyshev polynomials 
of the second kind.

\medskip
\noindent{{\sc 2001 Mathematics Subject Classification: Primary 05A05, 05A15; 
Secondary 30B70 42C05}}
\section*{\centering{\sc 1. Introduction}}

Let $\pi\in S_n$ and $\tau\in S_k$ be two permutations.
An {\it occurrence\/} of $\tau$ in $\pi$
is a subsequence $1\ls i_1<i_2<\dots<i_k\ls n$ such that $(\pi_{i_1},
\dots,\pi_{i_k})$ is order-isomorphic to $\tau$; in such a context $\tau$ is
usually called a {\it pattern\/}. We say that $\pi$ {\it avoids\/} $\tau$,
or is $\tau$-{\it avoiding\/}, if there is no occurrence of $\tau$ in $\pi$.
The set of all $\tau$-avoiding permutations in $S_n$ is denoted 
$S_n(\tau)$. \\

Pattern avoidance proved to be a useful language in a variety of seemingly 
unrelated problems, from stack sorting ~\cite{Kn}, ~\cite{Rt} to theory of 
Kazhdan-Lusztig polynomials ~\cite{Fb}, and singularities of Schubert 
varieties ~\cite{LS}, ~\cite{SCb}.
A natural generalization of single pattern avoidance is 
{\it subset avoidance\/}; that is, we say that $\pi\in S_n$ 
avoids a subset $T\subset S_k$ if $\pi$ avoids any $\tau\in T$. 
A complete study of subset avoidance for the case $k=3$ is
carried out in ~\cite{SS} (see also \cite{W,M1,M2}). \\

Several recent papers \cite{CW,RWZ,MV1,Kr,JR,MV2,MV3,MV4} deal with the case 
$\tau_1\in S_3$, $\tau_2\in S_k$ for various pairs $\tau_1,\tau_2$. Another
natural question is to study permutations avoiding $\tau_1$ and containing
$\tau_2$ exactly $t$ times. Such a problem for certain $\tau_1,\tau_2\in S_3$ 
and $t=1$ was investigated in \cite{R}, and for certain $\tau_1\in S_3$, 
$\tau_2\in S_k$ in \cite{CW,RWZ,MV1,Kr,JR,MV2,MV3,MV4}. The tools involved 
in these papers include continued fractions, Chebyshev polynomials, 
Dyck paths, and ordered trees.

\begin{definition}
{\it A finite continued fraction with $n$ steps} define as 
the following expression
	$$\frac{a_1}{\displaystyle b_1+\frac{a_2}{\displaystyle b_2+\ddots_{\displaystyle \frac{a_n}{b_n+a_{n+1}}}}}.$$
\end{definition}

There are many faces for applications of theory of continued fractions as an examples:  
Theory of functions, Approximation theory, Numerical analysis, and Restricted pattern.
As an application in restricted pattern is appear the continued fraction 
	$$\frac{1}{\displaystyle 1-\frac{x}{\displaystyle 1-\ddots_{\displaystyle \frac{x}{1-x}}}}.$$
in \cite{RWZ,CW}, and later than in \cite{MV1,Kr,JR,MV2,MV3,MV4}. Now we generalize this continued 
fraction by the following.

\begin{definition}
Let us denote the continued fraction with $k$ steps
	$$\frac{1}{\displaystyle 1-\frac{x}{\displaystyle 1-\ddots_{\displaystyle \frac{x}{1-xE}}}}$$
by $R_{k;E}(x)$ for any $k\geq 1$, and for $k=0$ we define $R_{0;E}(x)=E$. Also 
for simplicity we denote $R_{k;0}(x)$ by $R_k(x)$.
\end{definition}

Properties for $R_{k;E}$ is given by the following proposition.

\begin{proposition}
\label{REK}
Let $E$ any expression. Then
\begin{itemize}
\item[(i)]	For all $k\geq 1$; 
		$$R_{k;E}(x)=\frac{U_{k-1}\tx-\sqrt{x}E\cdot U_{k-2}\tx}{\sqrt{x}\left( U_k\tx-\sqrt{x}E\cdot U_{k-1}\tx\right) },$$
		where $U_k$ if the $k$th Chebyshev polynomials of the second kind;

\item[(ii)]	For all $k\geq 1$; 
		$$\prod\limits_{j=1}^k R_{j;E}(x)=\frac{1}{x^{\frac{k}{2}}\left[ U_k\tx -\sqrt{x}E\cdot U_{k-1}\tx \right] },$$
		where $U_k$ is the $k$th Chebychev polynomial of the second kind;

\item[(iii)]	$$\lim\limits_{k\rightarrow\infty} R_{k;E}(x)=\frac{1-\sqrt{1-4x}}{2x}.$$
\end{itemize}
\end{proposition}
\begin{proof}
$(i).$ For $k=1$ the proposition is trivial. By definitions 
	$$R_{k+1;E}(x)=\frac{1}{1-xR_{k;E}(x)},$$
and by us induction we yields 
	$$R_{k+1;E}(x)=\frac{\sqrt{x}\left( U_k\tx-\sqrt{x}E\cdot U_{k-1}\tx\right)}
		{ \sqrt{x}\left( U_k\tx-\sqrt{x}E\cdot U_{k-1}\tx\right) -x\left( U_{k-1}\tx-\sqrt{x}E\cdot U_{k-2}\tx \right)},$$
which means that 
	$$R_{k+1;E}(x)=\frac{U_k\tx-\sqrt{x}E\cdot U_{k-1}\tx}
		{ U_k\tx-\sqrt{x}E\cdot U_{k-1}\tx -\sqrt{x}\left( U_{k-1}\tx-\sqrt{x}E\cdot U_{k-2}\tx \right)}.$$
On the other hand, by definition of Chebyshev polynomials of the second kind have the following 
property $U_k(t) =2t U_{k-1}(t) - U_{k-2}(t)$, so 
	$$\sqrt{x}U_k\tx= U_{k-1}\tx -\sqrt{x}U_{k-2}\tx.$$
Hence the Proposition holds for $k+1$. \\

Again by us induction it is easy to see the second property, and the third property its yield 
immediately from \cite[Lemma 3.1]{MV1}.
\sof\end{proof}

\begin{example}
\label{exx}
By Proposition \ref{REK} $R_{k;0}=R_k(x)$, $R_{k,1}(x)=R_{k+1}(x)$, and 
$$R_{k;1+x}(x)=\frac{U_k\tx-xU_{k-2}\tx}{\sqrt{x}\left[ U_{k+1}\tx -xU_{k-1}\tx\right]}.$$
\end{example}

Now, for any three set of patterns $A$, $B$, and $C$ let us define $f_{A;B}^C(n)$ be the 
number of all $\alpha\in S_n(A)$ such that $\alpha$ containing every pattern in $B$ 
exactly once and containing every pattern in $C$ at least once. The corresponding 
generating function we denote by $F_{A;B}^C(x)$. For simplicity, we write $F_A(x)$, 
$F_{A;B}(x)$ and $G_B(x)$ when $B=C=\emptyset$, $C=\emptyset$, and $A=C=\emptyset$ 
respectively. \\

The paper is organized as the following. In section $2$, we find a recurrence in terms of generating functions for $F_{A;B}^C(x)$,  
and we prove $F_{A;B}^C(x)$ for all $A,B,C$ such that $A\cup B\neq\emptyset$ 
is a rational function. In section $3$ we present an examples of the main results, which 
are present the relation between the restricted patterns and continued fractions. 
\section*{\centering{\sc 2. Main results}}
Consider an arbitrary pattern $\tau=(\tau_1,\dots,\tau_k)\in S_k(132)$. 
Recall that $\tau_i$ is said to be a {\it right-to-left maximum\/} if 
$\tau_i>\tau_j$ for any $j>i$. Let $m_0=k,m_1,\dots,m_r$
be the right-to-left maxima of $\tau$ written from left to right. Then $\tau$
can be represented as
		$$\tau=(\tau^0,m_0,\tau^1,m_1,\dots,\tau^r,m_r),$$
where each of $\tau^i$ may be possibly empty, and all the entries of $\tau^i$
are greater than $m_{i+1}$ and all the entries of $\tau^{i+1}$. 
This representation is called
the {\it canonical decomposition\/} of $\tau$. Given the canonical
decomposition, we define the $i$th {\it prefix\/} of $\tau$ by 
$\pi^i(\tau)=(\tau^0,m_0,\dots,\tau^i,m_i)$ for $1\leq i\ls eq$ and $\pi^0(\tau)=\tau^0$,
$\pi^{-1}(\tau)=\emptyset$.
Besides, the $i$th {\it suffix\/} of $\tau$ is defined by
$\sigma^i(\tau)=(\tau^i,m_i,\dots,\tau^r,m_r)$ for $0\leq i\leq r$ and
$\sigma^{r+1}(\tau)=\emptyset$. Strictly speaking, prefixes and suffices 
themselves are not patterns, since they are not permutations (except for
$\pi^r(\tau)=\sigma^0(\tau)=\tau$). However, 
any prefix or suffix is order-isomorphic to a unique permutation, and in what 
follows we do not distinguish between a prefix (or suffix) and the 
corresponding permutation.
Now, let us find $f_A^C(n)$ in terms of $f_T(n)$ by the following lemma.

\begin{lemma}
\label{s11}
Let $A$ set of patterns, and let $\{\beta^{(i)}\}_{i=1}^m$ be sequence 
of patterns. Then 
	$$f_A^{\beta^{(1)},\dots,\beta^{(m)}}(n)=
	\sum\limits_{j=0}^m \left( (-1)^j \sum\limits_{1\leq i_1<i_2<\dots<i_j\leq m} 
	f_{A,\beta^{(i_1)},\dots,\beta^{(i_j)}}(n) \right).$$
\end{lemma}
\begin{proof}
By definitions 
	$$f_A^{\beta^{(1)}}(n)+f_{A,\beta^{(1)}}(n)=f_A(n),$$
which means this statement holds for $m=1$. So more generally, 
	$$f_A^{\beta^{(1)},\dots,\beta^{(m+1)}}(n)=
	\sum\limits_{j=0}^m \left( (-1)^j \sum\limits_{1\leq i_1<i_2<\dots<i_j\leq m} 
	f_{A,\beta^{(i_1)},\dots,\beta^{(i_j)}}^{\beta^{(m+1)}}(n) \right),$$
by use induction and the same argument in the case $m=1$, the theorem holds.
\sof\end{proof}

Immediately, we can represent this result (Lemma \ref{s11}) by another way, as the 
following.

\begin{theorem}
\label{c11}
Let $\{\alpha^{(i)}\}_{i=1}^m$, $\{\beta^{(i)}\}_{i=1}^m$ be two 
sequences of patterns such that $\alpha^{(i)}$ contains $\beta^{(i)}$ 
for all $i=1,2,\dots,m$, and let $A$ set of patterns. Then 
	$$f_{A,\alpha^{(1)},\dots,\alpha^{(m)}}^{\beta^{(1)},\dots,\beta^{(m)}}(n)=
	\sum\limits_{j=0}^m \left( (-1)^j \sum\limits_{1\leq i_1<i_2<\dots<i_j\leq m} 
	f_{A_{i_1,\dots,i_j}}(n) \right),$$
where $A_{i_1,\dots,i_j}=A \bigcup\limits_{d=1}^m\{\alpha^{(d)}\}
	\backslash\{\alpha^{(i_1)},\dots,\alpha^{(i_j)}\}
	\bigcup\limits_{d=1}^j\{\beta^{(i_d)}\}.$
\sof\end{theorem}

Now we present the main two results of this paper in the following two subsections. 
The first result, we find a recurrence to calculate the generating function $F_T(x)$ 
where $T$ set of pattern. The second result, we find an another recurrence to calculate 
the generating function $F_{A;B}(x)$, which is a generalization of the first result.

\subsection*{The generating function $F_A(x)$.} 
In the current subsection we find a recurrence to calculate $F_A(x)$. This calculation 
immediately by induction, by uses Theorem \ref{c11}, and us result in \cite[Th 2.1]{MV3}.

\begin{theorem}
\label{avoidsetp}
Let $\tau^{[i]}=(\tau^{i,0},d_{i,0},\tau^{i,1},d_{i,1},\dots,\tau^{i,m_i},d_{i,m_i})\in S_{d_{i,0}}(132)$ for $i=1,2,\dots,p$ 
such no there two patterns one contain the another. Then 
$$F_{\tau^{[1]},\dots,\tau^{[p]}}(x)=
	1+x\sum\limits_{j_1=0}^{m_1}\sum\limits_{j_2=0}^{m_2}\dots\sum\limits_{j_p=0}^{m_p} 
		F_{\{\pi_{\tau^{[1]}}^{(j_1)},\dots,\pi_{\tau^{[p]}}^{(j_p)}\}}^{\pi_{\tau^{[1]}}^{(j_1-1)},\dots,\pi_{\tau^{[p]}}^{(j_p-1)}}(x)
		F_{\{\theta_{\tau^{[1]}}^{(j_1)},\dots,\theta_{\tau^{[p]}}^{(j_p)}\}}(x),$$
where 
	$$F_{\{\pi_{\tau^{[1]}}^{(j_1)},\dots,\pi_{\tau^{[p]}}^{(j_p)}\}}^{\pi_{\tau^{[1]}}^{(j_1-1)},\dots,\pi_{\tau^{[p]}}^{(j_p-1)} }(x)
	=\sum\limits_{j=0}^p (-1)^j \sum\limits_{1\leq i_1<\dots<i_j\leq p} F_{A_{i_1,\dots,i_j}^{j_1,\dots,j_p}}(x),$$
such that $A_{i_1,\dots,i_j}=\cup_{q\neq i_1,\dots,i_j} \{\pi_{\tau^{[q]}}^{(j_q)}\} 
				     \bigcup \cup_{q= i_1,\dots,i_j} \{\pi_{\tau^{[q]}}^{(j_q-1)}\}.$
\end{theorem}

\subsection*{The generating function $F_{A;B}(x)$.} 
Here, we find a recurrence to calculate $F_{A;B}(x)$. This calculation 
immediately by induction, by uses Theorem \ref{c11}, Theorem \ref{avoidsetp}, 
and us result in \cite[Th 3.1]{MV3}.

\begin{theorem} 
\label{avoidcontain}
Let $A,B$ any disjoint sets of pattern. Then 
\begin{itemize}
\item[(i)]
$F_{132,A;B}(x)$ is a rational function satisfying the relation
$$
F_{A;B}(x)=	\left\{ 
		\begin{array}{ll}
		x\sum\limits_{a_i=0,1,\dots,r_{\tau_i}}
			\sum\limits_{b_j=0,1,\dots,1+r_{\gamma_j}}
			F_{A_1;B_1}^C(x)F_{A_2;B_2}(x),& \ B\neq\emptyset \\
		1+x\sum\limits_{a_i=0,1,\dots,r_{\tau_i}}
			\sum\limits_{b_j=0,1,\dots,1+r_{\gamma_j}} 
			F_{A_1;B_1}^C(x)F_{A_2;B_2}(x),& \ B=\emptyset 
		\end{array}
		\right.
$$
where $A=\{ \tau_i | i=1,2,\dots,a\}$, $B=\{ \gamma_j | j=1,2,\dots,b\}$, and 
$$\begin{array}{l}
	C=\{\pi^{a_i-1}(\tau_i)| i=1,2,\dots,a\}, \\
	A_1=\{\pi^{a_i}(\tau_i)| i=1,2,\dots,a\}\cup \{ \pi^{b_j}(\gamma_j) |j=1,2,\dots,b\}, \\
	B_1=\{\pi^{b_j-1}(\gamma_j)| j=1,2,\dots,b\}, \\
	A_2=\{\sigma^{a_i}(\tau_j)| i=1,2,\dots,a\}\cup\{\sigma^{b_j-1}(\gamma_j)| j=1,2,\dots,b\},\\
	B_2=\{\sigma^{b_j}(\gamma_j)| j=1,2,\dots,b\}.
\end{array}$$

\item[(ii)] for two sequences of patterns $\{\alpha^{(i)}\}_{i=1}^m$, $\{\beta^{(i)}\}_{i=1}^m$ 
such that $\alpha^{(i)}$ contains $\beta^{(i)}$ 
for all $i=1,2,\dots,m$, and for any set of patterns $T$, 
	$$F_{T,\alpha^{(1)},\dots,\alpha^{(m)};B}^{\beta^{(1)},\dots,\beta^{(m)}}(x)=
	\sum\limits_{j=0}^m \left( (-1)^j \sum\limits_{1\leq i_1<i_2<\dots<i_j\leq m} 
	F_{A_{i_1,\dots,i_j};B}(x) \right),$$
where $A_{i_1,\dots,i_j}=T\bigcup\limits_{d=1}^m\{\alpha^{(d)}\}
	\backslash\{\alpha^{(i_1)},\dots,\alpha^{(i_j)}\}
	\bigcup\limits_{d=1}^j\{\beta^{(i_d)}\}.$
\end{itemize}
\end{theorem}
\section*{\centering{\sc 3. Examples and continued fractions}}

Though elementary, Theorem \ref{avoidsetp} enables us to derive easily various 
known and new results for a fixed a set of patterns. 

\begin{example} 
\label{exa1}
{\rm (see \cite{G})} An a numerical case, by Theorem \ref{avoidsetp} we yields 
$$\begin{array}{l}
F_{\{2341,3241\}}(x)= 1+xF_{\{23,32\}}(x)F_{\{2341,3241\}}(x)+\\
\ \ \ +x(F_{\{23,3241}(x)-F_{\{23,32\}}(x))F_{\{2341,1\}}(x)+\\
\ \ \ +x(F_{\{2341,32\}}(x)-F_{\{23,32\}}(x))F_{\{1,3241\{}(x)+\\
\ \ \ +x(F_{\{2341,3241\}}(x)-F_{\{23,3241\}}(x)-F_{\{2341,32\}}(x)+F_{\{23,32\}}(x))F_{\{1,1\}}(x).
\end{array}$$
On the other hand, by definition it is easy to see $F_{\{23,32\}}(x)=1+x$, $F_{\{\tau,1\}}(x)=1$,  
and $F_{\{23,3241\}}(x)=F_{\{2341,32\}}=\frac1{1-x}$. Hence,  
it is easy to get $F_{\{2341,3241\}}(x)=\frac{1-x-x^2}{1-2x-x^2}$.
\end{example}

As a corollary of Theorem \ref{avoidsetp} we obtain the following.

\begin{corollary}
\label{endk}
Let $T'$ set of pattern, and let 
$T=\{(\tau_1,\dots,\tau_{k-1},k)| (\tau_1,\dots,\tau_{k-1})\in T'\}$. Then 
	$$F_T(x)=\frac{1}{1-xF_{T'}(x)}.$$
\end{corollary}

\begin{example} 
\label{exa2} 
{\rm (see \cite[Pr. 15]{SS}, and \cite[Sec 4.1]{G})}
An another numerical example, by Corollary \ref{endk} we yields
	$$F_{\{123,213\}}(x)=\frac{1}{1-xF_{\{21,12\}}(x)},$$
and by definitions we get the result \cite[Pr. 15]{SS}, which is 
	$$F_{\{123,213\}}(x)=\frac{1}{1-x-x^2}.$$
In the same way, by use Corollary \ref{endk} twice we yields 
	$$F_{\{1234,2134\}}(x)=\frac{1-x-x^2}{1-2x-x^2},$$
which is result \cite[Sec. 4.1]{G}. 
\end{example}

Now let us generalize the above example. First let define a special set of patterns.
 
\begin{definition}
For any $k\geq l\geq 1$, let $U_l^k$ be the set of all permutations 
$\tau\in S_k$ such that $(\tau_{l+1},\tau_{l+2},\dots,\tau_k)=(l+1,l+2,\dots,k)$.
Clearly $|U_l^k|=l!$.
\end{definition}

By \cite{Kn} and definitions  
		$$F_{U_l^l}(x)=\sum_{j=0}^{l-1} c_jx^j,$$
where $c_j$ is the $j$th Catalan number. Hence, consequentially to Example \ref{exa2} 
we yields similarly the following.

\begin{corollary}
\label{ce1}
Let $k\geq l\geq 1$; then 
	$$F_{U_l^k}(x)=R_{k-l;E(x)}(x),$$
where $E(x)=\sum\limits_{j=0}^{l-1} c_jx^j$, and $c_j$ is the $j$th Catalan number.
\end{corollary}

\begin{example} {\rm(see \cite{CW,MV1,Kr})}
For $l=1$, by Corollary \ref{ce1} we yields 
	$$F_{U_1^k}(x)=F_{12\dots k}(x)=R_{k;0}(x)=R_k(x).$$
\end{example}

Now we present another direction to use continued fractions.

\begin{corollary}
\label{ce2}
Let $k>l\geq 1$. For any $\tau\in U_l^k$, 
	$$F_{U_l^k\backslash\{\tau\};\tau}(x)=\frac{x^l}{\left( U_{k-l}\tx -\sqrt{x}E(x)U_{k-1-l}\tx\right)^2}.$$
where $E(x)=\sum\limits_{j=0}^{l-1} c_jx^j$, and $c_j$ is the $j$th Catalan number.
\end{corollary}
\begin{proof}
Let us fix $\tau\in U_l^k$ such that $(\tau=\tau',k)$, 
and let us denote $U_l^k\backslash\{\tau\}$ by $M_l^k$.
Immediately by Theorem \ref{avoidcontain} we yields
$$F_{M_l^k;\tau}(x)=xF_{M_l^{k-1};\tau'}(x)F_{U_l^k}(x)+xF_{U_l^{k-1}}(x)F_{M_l^k;\tau}(x),$$
which means by Corollary \ref{ce1}
	$$F_{M_l^k;\tau}(x)=\frac{xR_{k-l;E}(x)}{1-xR_{k-1-l;E}(x)} F_{M_l^{k-1};\tau'}(x),$$
where $E=\sum\limits_{j=0}^{l-1} c_jx^j$, and $c_j$ is the $j$th Catalan number.
Since $\frac{1}{1-xR_{m-1;E}(x)}=R_{m;E}(x)$ we obtain 
 	$$F_{M_l^k;\tau}(x)=xR_{k-l;E}^2(x) F_{M_l^{k-1};\tau'}(x).$$
By use induction we yields 
 	$$F_{M_l^k;\tau}(x)=x^{k-l}F_{M_l^l;\beta}(x) \prod_{j=1}^{k-l} R_{j;E}^2(x),$$
where $\tau=(\beta,l+1,l+2,\dots,k)$. On the other hand, by definitions $F_{M_l^l;\beta}(x)=x^l$, so  
 	$$F_{M_l^k;\tau}(x)=x^k \prod_{j=1}^{k-l} R_{j;E}^2(x).$$
Hence by Proposition \ref{REK} the corollary holds.
\sof\end{proof}

\begin{example} {\rm (see \cite{MV1,Kr})} 
\label{exc1}
For either $l=1$, $E=1$; or $l=0$, $E=0$, we yields from 
Corollary \ref{ce2} for all $k\geq 1$ 
	$$G_{\{12\dots k\}}(x)=F_{\{\emptyset;12\dots k\}}(x)=\frac{1}{U_k^2\tx }.$$
For $l=2$, $E=1+x$ we obtain for all $k\geq 3$ 
	$$F_{\{123\dots k;213\dots k\}}(x)=F_{\{213\dots k;123\dots k\}}(x)=\frac{x}{\left( U_{k-1}\tx -xU_{k-3}\tx\right)^2}.$$
\end{example}

\begin{example}
Now we complete all the calculation for either containing exactly once, or avoiding 
a two patterns from $U_2^k=\{123\dots k, 213\dots k\}$. By Corollary \ref{ce1} and 
Example \ref{exc1} its left to find $G_{U_2^k}(x)$. Let $k\geq 3$; 
by Theorem \ref{avoidcontain} 
{\small 
$$\begin{array}{ll} 
G_{U_2^k}(x)	&=xF_{U_2^{k-1}}(x)G_{U_2^k}(x)+xF_{123\dots k-1;213\dots k-1}(x) F_{213\dots k; 123\dots k}(x)\\
		&+xF_{213\dots k-1;123\dots k-1}(x) F_{123\dots k; 213\dots k}(x)+xG_{U_2^{k-1}}(x)F_{U_2^k}(x),
\end{array}$$}
so by Corollary \ref{ce1}, Example \ref{exc1}, and  by use Proposition \ref{REK} we yields 
$$G_{U_2^k}(x)	=\frac{2x^2\sqrt{x}}{W_{k;2}(x)W_{k;1}^2(x)} 
		+xR_{k-2;1+x}^2(x)G_{U_2^{k-1}}(x),$$ 
where $W_{k;j}(x)=U_{k-j}\tx -xU_{k-2-j}\tx$. 
Besides $G_{U_2^k}(x)=0$ for all $k=0,1,2,3$, hence 
$$G_{U_2^k}(x)	=\frac{2x^2\sqrt{x}}{W_{k;1}^2(x)} \sum\limits_{j=3}^{k-2} 
		\frac{1}{W_{k;j-1}(x)W_{k;j}(x)}.$$
\end{example}

\end{document}